\documentclass[11pt,bezier]{article}
\pdfoutput=1
\usepackage{array}
\usepackage{multirow}
\usepackage{makecell}
\usepackage{graphicx}
\usepackage{booktabs}
\usepackage{colortbl}
\usepackage{amsmath}
\usepackage{amsfonts,amsthm,amssymb}
\usepackage{exscale}
\usepackage{relsize}
\usepackage{graphicx}
\usepackage{epstopdf}
\usepackage{color}
\usepackage{ifpdf}
\usepackage{caption}
\usepackage{float}
\usepackage{enumitem}
\usepackage{hyperref}
\usepackage{amssymb}
\usepackage{amsmath}
\usepackage{amsthm}
\newcounter{mycounter}
\theoremstyle{plain}
\newtheorem{lemma}[mycounter]{Lemma}
\newtheorem{theorem}[mycounter]{Theorem}
\newtheorem{observation}[mycounter]{Observation}
\newtheorem{proposition}[mycounter]{Proposition}
\newtheorem{corollary}[mycounter]{Corollary}
\newtheorem{conjecture}[mycounter]{Conjecture}
\newtheorem{definition}[mycounter]{Definition}
\makeatletter
\@addtoreset{mycounter}{section}
\makeatother

\usepackage{graphicx}
\makeatletter
\newcommand{\Rmnum}[1]{\expandafter\@slowromancap\romannumeral #1@}
\makeatother
\theoremstyle{definition}

\usepackage[justification=centering]{caption}
\allowdisplaybreaks[3]
\begin{document}
%
%

\title{On the Steiner $k$-diameter and Steiner ($k,k^{\prime}$)-radius of trees\footnote{The research is supported by National Natural Science Foundation of China (12261086) and Natural Science Foundation of Xinjiang Uygur Autonomous Region (2025D01E02).}}
\author{Qingnan Zhang, Yingzhi Tian\footnote{Corresponding author. E-mail: 18738907528@163.com (Q. Zhang); tianyzhxj@163.com (Y. Tian).} \\
{\small College of Mathematics and System Sciences, Xinjiang
University, Urumqi, Xinjiang 830046, PR China}}

\date{}
\maketitle

\begin{sloppypar}

\noindent{\bf Abstract } Given a connected graph $G=(V,E)$ and a $k$-set $S\subseteq V(G)$, the $Steiner$ $distance$ $d_{G}(S)$ of $S$ is defined as the size of a minimum tree including $S$ in $G$.
The $Steiner$ $k$-$eccentricity$ of a vertex $v$ in $G$ is the maximum value of $d_G(S)$ over all $S\subseteq V(G)$ with $|S|=k$ and $v\in S$.
The minimum Steiner $k$-eccentricity over all vertices, denoted by $Sr_k(G)$, is called the $Steiner$ $k$-$radius$ of $G$ and the maximum Steiner $k$-eccentricity over all vertices, denoted by $Sd_k(G)$, is its $Steiner$ $k$-$diameter$.
The $Steiner$ $(k,k^{\prime})$-$eccentricity$ of a $k^{\prime}$-subset $S^{\prime}$ of $V(G)$, which is an extension of the Steiner $k$-eccentricity of a vertex $v$, is defined as the maximum Steiner distance over all $k$-subsets of $V(G)$ containing $S^{\prime}$.
The minimum Steiner $(k,k^{\prime})$-eccentricity among all $k^{\prime}$-subsets of $V(G)$, denoted by $Sr_{k,k^{\prime}}(G)$, is called the $Steiner$ $(k,k^{\prime})$-$radius$ of $G$.
In 1989, Chartrand, Oellermann, Tian and Zou showed that for any $k\geq3$, $Sd_k(T)\leq \frac{k}{k-1}Sr_k(T)$ for any tree $T$.
In this paper, we generalize this result and show that $Sd_k(T)\leq \frac{k}{k-k^{\prime}}Sr_{k,k^{\prime}}(T)$ for any $k\geq3$, $k>k^{\prime}\geq1$.
Furthermore, for $k^{\prime}=2$ and $k^{\prime}=3$, we obtain a tight upper bound of the Steiner $k$-diameter by the Steiner $(k,k^{\prime})$-radius for all trees.


\section{Introduction}
~~~~All graphs considered in this paper are finite, simple and connected.
Let $G$ be a graph with vertex set $V(G)$ and edge set $E(G)$. The $order$ of $G$ is $|V(G)|$, and its $size$ is $|E(G)|$.
The $distance$ between two vertices $u, v\in V(G)$ is defined as the length of the shortest path connecting them, denoted by $d_{G}(u,v)$. If $H$ is a subgraph of $G$ and $v\in V(G)$, then the distance from $v$ to $H$, is defined as $d_G(v,H)=\min\{d_G(v,u): u\in V(H)\}$.
The $eccentricity$ $\varepsilon(v)$ of a vertex $v$ in $G$ is defined as
$\varepsilon(v)=\max\{d(v, u): u\in V(G)\}$.
The $radius$ and the $diameter$ of the graph $G$ are $rad(G)=\min\{\varepsilon(v): v\in V(G)\}$ and $diam(G)=\max\{\varepsilon(v): v\in V(G)\}$, respectively.
These last two parameters satisfy the inequalities $rad(G)\leq diam(G)\leq 2rad(G)$.
A $diametrical$ $path$ is the shortest path between two vertices with distance $diam(G)$.
In a graph with radius $rad(G)$, a $central$ $vertex$ is identified as a vertex with eccentricity equal to $rad(G)$. The set of central vertices is denoted by $C(G)$. It is evident that any tree contains either one or two central vertices. In a graph with diameter $diam(G)$, a $peripheral$ $vertex$ is defined as a vertex with eccentricity equal to $diam(G)$.

The Steiner distance was first introduced in 1989 by Chartrand, Oellermann, Tian and Zou \cite{Steiner distance}. The Steiner distance is an extension of the classical graph distance. Correspondingly, some parameters related to the classical graph distance can be extended to the corresponding parameters involved in the Steiner distance.
The $Steiner$ $distance$ $d_{G}(S)$ among the vertices of $S$ (or simply the distance of $S$) is defined as the minimum size among all connected subgraphs whose vertex sets contain $S$. Notably, if $H$ is a connected subgraph of $G$ such that $S\subseteq V(H)$ and $|E(H)|=d_{G}(S)$, then $H$ must be a tree, that is,
$$d_{G}(S)=\min\{|E(T)|: S\subseteq V(T)\},$$
where $T$ is a subtree of $G$.
If $T$ satisfies $S\subseteq V(T)$ and $d_{G}(S)=|E(T)|$, then $T$ is an $S$-Steiner tree of $G$. We will denote this tree as $T_{S}$.
The $Steiner$ $k$-$eccentricity$ $\varepsilon_{k}(v;G)$ of a vertex $v$ in $G$ is defined by $\max\{d_G(S): v\in S\subseteq V(G), |S|=k\}$, where $k\geq2$. The Steiner $k$-radius of $G$ is $Sr_k(G)=\min\{\varepsilon_k(v): v\in V(G)\}$, and the $Steiner$ $k$-$diameter$ of $G$ is $Sd_{k}(G)=\max\{\varepsilon_k(v): v\in V(G)\}$.
In 1989, Chartrand, Oellermann, Tian and Zou \cite{Steiner distance} derived the sharp upper bound of the Steiner $k$-diameter by the Steiner $k$-radius of trees. In 1990, Henning, Oellermann and Swart \cite{diameter radius1} derived the sharp upper bound of the Steiner $k$-diameter by the Steiner $k$-radius of general connected graphs for $k=3, 4$.
In 2024, Reiswig \cite{diameter radius2} obtained the upper bound of the Steiner $k$-diameter by the Steiner $k$-radius of general connected graphs when $k\geq5$ and verified that this bound is tight.

The upper bounds on the Steiner $k$-radius in terms of the order and the minimum degree for all graphs and for graphs that contain no triangles were provided in \cite{steiner radius}. The upper bounds on the Steiner $k$-diameter in terms of the order, the minimum degree, and the girth of a graph were developed in \cite{steiner diameter bound1,steiner diameter bound2}.
The bounds on the Steiner $k$-diameter for specific graph classes, such as maximal planar graphs, hypercubes, product graphs, were established in \cite{connected,hypercube,product}.
Furthermore, the bounds on the size of graphs with given Steiner $k$-diameter and maximum degree were given in \cite{Steiner diameter2,Steiner diameter4}.
The graphs with $Sd_{n-k}(G)=l$ for $n-k-1\leq l\leq n-1$, $Sd_3(G)=2,3,n-1$ and $Sd_4(G)=3,4,n-1$ were characterized in \cite{Steiner diameter3,Steiner diameter,Steiner diameter1}.
The relationship between the Steiner $k$-diameter of a graph and the Steiner $k$-diameter of the corresponding line graph was described in \cite{line}.
Additionally, for the study of some parameters related to the Steiner distance, such as the Steiner $k$-eccentricity, the Steiner centers and the Steiner medians, see \cite{Steiner 3-eccentricity4,Steiner 3-eccentricity3,Steiner 3-eccentricity1,center}.

The Steiner $(k,k^{\prime})$-eccentricity, an extension of the Steiner $k$-eccentricity, was first introduced in 2024 by Li, Zeng, Pan and Li \cite{average Steiner eccentricity32}. The $Steiner$ $(k, k^{\prime})$-$eccentricity$ $\varepsilon_{k,k^{\prime}}(S^{\prime};G)$ of a $k^{\prime}$-subset $S^{\prime}$ in $G$ is defined as the maximum Steiner distance over $k$-subsets of $V(G)$ that include $S^{\prime}$, that is,
$$\varepsilon_{k,k^{\prime}}(S^{\prime}; G)=\max\{d_{G}(S): S^{\prime}\subseteq S\subseteq V(G), |S|=k\},$$
where $k\geq2$ and $k\geq k^{\prime}\geq1$.
Consequently, the Steiner $k$-eccentricity $\varepsilon_k(v;G)$ of $v\in V(G)$ is exactly $\varepsilon_{k,1}(v;G)$.
A $k$-subset $S$ of $V(G)$ is a Steiner $k$-ecc $S^{\prime}$-set if it satisfies $S^{\prime}\subseteq S$ and $d_{G}(S)=\varepsilon_{k,k^{\prime}}(S^{\prime};G)$, where $|S^{\prime}|=k^{\prime}$. The corresponding tree $T$ that realizes $\varepsilon_{k,k^{\prime}}(S^{\prime};G)$ is called a Steiner $k$-ecc $S^{\prime}$-tree. It is noteworthy that the Steiner $k$-ecc $S^{\prime}$-set may not be unique, and each such set $S^{\prime}$ can potentially yield more than one Steiner $k$-ecc $S^{\prime}$-tree.
The $Steiner$ $(k,k^{\prime})$-$radius$ of $G$ is $Sr_{k,k^{\prime}}(G)=\min\{\varepsilon_{k,k^{\prime}}(S^{\prime};G): S^{\prime}\subseteq V(G), |S^{\prime}|=k^{\prime}\}$, where $k^{\prime}\geq1$.
Note that the Steiner $k$-radius of $G$ is just $Sr_{k,1}(G)$. In this paper, we will generalize the results in \cite{Steiner distance} and give the upper bound of the Steiner $k$-diameter by the Steiner $(k,k^{\prime})$-radius of trees for $k\geq3$ and $k>k^{\prime}\geq1$.

The structure of the paper is as follows. In Section 2, we present fundamental results regarding the Steiner $k$-diameter and the Steiner $(k,k^{\prime})$-radius and also provide several existing theorems related to the Steiner $k$-diameter and the Steiner $k$-radius.
In Section 3, for $k\geq3$, $k>k^{\prime}\geq 1$, we obtain an upper bound of the Steiner $k$-diameter of trees by its Steiner $(k,k^{\prime})$-radius. In Section 4, we derive a sharp upper bound of the Steiner $k$-diameter by the Steiner $(k,k^{\prime})$-radius of trees for $k^{\prime}=2, 3$.

\section{Preliminary}

~~~~A vertex with degree 1 in $G$ is termed a $pendant$ $vertex$. The set of all pendant vertices is denoted as $L(G)$. Every tree of order greater than one contains at least two pendant vertices.
A vertex with degree at least 2 is defined as an $internal$ $vertex$ of $G$, and the set of internal vertices is denoted as $NL(G)$. Furthermore, a vertex with degree at least 3 is called a $branching$ $vertex$ of $G$.
The edge incident with a pendant vertex is referred to as a $pendant$ $edge$.
A path $P$ is defined as a $pendant$ $path$ in $G$ if one end vertex of $P$ is a pendant vertex, the other end vertex is a branching vertex, and the degree of each internal vertex of $P$ is exactly 2 in $G$. A $starlike$ $tree$ is a tree with exactly one branching vertex.

We first introduce some basic results of the Steiner $k$-diameter and the Steiner ($k,k^{\prime}$)-radius.

\begin{observation}
Let $n$, $k$ be two integers with $n\geq k\geq2$ and let $G$ be a connected graph of order $n$.

(1) If $H$ is a connected spanning subgraph of $G$, then $Sd_{k}(G)\leq Sd_{k}(H)$.

(2) $Sd_{k}(G)\leq Sd_{k+1}(G)$.
\end{observation}

\begin{observation}(\cite{average Steiner eccentricity32})
Let $n$, $k$ be two integers with $n\geq k\geq2$ and $G$ be a connected graph of order $n$. The Steiner $(k,k^{\prime})$-radius satisfies
$$Sr_{k,1}(G)\geq Sr_{k,2}(G)\geq \cdots \geq Sr_{k,k}(G).$$
\end{observation}

\begin{lemma}(\cite{Steiner distance})\label{lem2.4}
Let $k\geq2$ be an integer, and $T$ be a tree with at least $k$ pendant vertices. If $S$ is a $k$-subset of $V(T)$ and $Sd_{k}(T)=d_{T}(S)$, then $S\subseteq L(T)$.
\end{lemma}

\begin{lemma}(\cite{rainbow})\label{lem2.5}
Let $n$, $k$ be two integers with $n\geq k\geq2$, and $G$ be a connected graph of order $n$.
Then $$k-1\leq Sd_k(G)\leq n-1,$$
the left equality holds if $G\cong K_n$ and the right equality holds if $G\cong P_n$.
\end{lemma}

\begin{lemma}(\cite{Steiner centers})\label{lem2}
Let $k\geq2$ be an integer, and $T$ be a tree with $k$ pendant vertices. For every integer $p$ with $2\leq p\leq k$, there exists a $p$-subset $S_p$ of $V(T)$ satisfying $d_{T}(S_p)=Sd_p(T)$   such that $S_2\subset S_3\subset\cdots \subset S_k$.
\end{lemma}

\begin{corollary}\label{cor2}
Let $k\geq2$ be an integer, and $T$ be a tree with $k$ pendant vertices. Then $Sd_{k}(T)\leq d+(k-2)\lfloor\frac{d}{2}\rfloor$, where $d$ is the diameter of $T$.
\end{corollary}
\noindent{\bf Proof.} Let $S=S_k$, where $S_k$ is the same as this in Lemma \ref{lem2}. Then $d_T(S)=Sd_k(T)$ and $T_S$ contains a diametrical path in $T$. Therefore, we can derive that $Sd_k(T)=d_T(S)\leq d+(k-2)\lfloor\frac{d}{2}\rfloor$. $\Box$

\begin{lemma}(\cite{steinertree})\label{steinertree}
Let $S^{\prime}$ be a $k^{\prime}$-subset in a tree $T$. If $|L(T)\setminus S^{\prime}|\geq k-k^{\prime}>0$, then given a Steiner $(k-1)$-ecc $S^{\prime}$-tree $T^{k-1}$,
there is a Steiner $k$-ecc $S^{\prime}$-tree $T^{k}$
such that $T^{k}$ contains this $T^{k-1}$.
\end{lemma}

\begin{proposition}(\cite{Steiner diameter,Steiner diameter4})
Let $n$, $k$ be two integers with $n\geq k\geq2$.

(1) For a complete graph $K_n$, $Sd_k(K_n)=k-1$.

(2) For a path $P_n$, $Sd_k(P_n)=n-1$.

(3) For a complete $r$-partite graph $K_{n_1,n_2,\ldots,n_r}(n_1\leq n_2\leq\cdots\leq n_r)$,
\[
Sd_k(K_{n_1,n_2,\ldots,n_r})=\left\{
\begin{array}{ll}
k, & \text{if } n_r\geq k; \\
k-1, & \text{if } n_r\leq k-1.
\end{array}
\right.
\]
\end{proposition}

\begin{proposition}
Let $n$, $k$ and $k^{\prime}$ be integers with $n\geq k\geq 2$ and $k\geq k^{\prime}\geq 1$.

(1) For a complete graph $K_n$, $Sr_{k,k^{\prime}}(K_n)=k-1$.

(2) For a path $P_n$,
\[Sr_{k,k^{\prime}}(P_n)=\left\{
\begin{array}{ll}
k-1, & \text{if } k-k^{\prime}=0 ;\\
\lceil \dfrac{n+k^{\prime}-2}{2}\rceil, & \text{if } k-k^{\prime}=1 ;\\
n-1, & \text{if } k-k^{\prime}\geq2.
\end{array}
\right.\]

(3) For a complete $r$-partite graph $K_{n_1,n_2,\ldots,n_r}(n_1\leq n_2\leq\cdots\leq n_r)$.

If $k^{\prime}=1$, then  \(Sr_{k,k^{\prime}}(K_{n_1,n_2,\ldots,n_r})=\left\{
\begin{array}{ll}
k, & \text{if } n_1\geq k ;\\
k-1, & \text{if } n_1\leq k-1.
\end{array}
\right.\)

If $k^{\prime}\geq2$, then $Sr_{k,k^{\prime}}(K_{n_1,n_2,\ldots,n_r})=k-1$.

\end{proposition}
\noindent{\bf Proof.} (1) For any $k$-subset $S$ of $V(K_n)$, we have $d_{K_n}(S)=k-1$. Hence, $Sr_{k,k^{\prime}}(K_n)=k-1$.

(2) Let $V(P_n)=\{v_0, v_1, \ldots, v_{n-1}\}$. If $k-k^{\prime}=0$, then $Sr_{k,k^{\prime}}(P_n)=Sr_{k,k}(P_n)$. Let $S^{\prime}=\{v_0,v_1,\ldots,v_{k-1}\}$. Then it follows that $\varepsilon_{k,k}(S^{\prime};P_n)=k-1$. Hence, $Sr_{k,k}(P_n)\leq k-1$. Moreover, since $Sr_{k,k}(P_n)\geq k-1$, we conclude that $Sr_{k,k}(P_n)=k-1$.

If $k-k^{\prime}=1$ and $n$ is odd, then there is only one central vertex $v_{(n-1)/2}$ in $P_{n}$. Let $S^{\prime}=\{v_{(n-1)/2-\lceil(k^{\prime}-1)/2\rceil}, \ldots, v_{(n-1)/2-1}, v_{(n-1)/2}, v_{(n-1)/2+1}, \ldots, v_{(n-1)/2+\lfloor(k^{\prime}-1)/2\rfloor}\}$ and $S=S^{\prime}\cup \{v_{n-1}\}$, we have  $\varepsilon_{k,k^{\prime}}(S^{\prime};P_n)=d_{P_n}(S)=\lceil\frac{k^{\prime}-1}{2}\rceil+\frac{n-1}{2}$.
Hence, $Sr_{k,k^{\prime}}(P_n)\leq \lceil\frac{k^{\prime}-1}{2}\rceil+\frac{n-1}{2}$.
In the following, we will show that $\varepsilon_{k,k^{\prime}}(S^{\prime \prime};P_n)\geq \lceil\frac{k^{\prime}-1}{2}\rceil+\frac{n-1}{2}$ for any $k^{\prime}$-subset $S^{\prime \prime}$.
Let $S^{\prime \prime}$ be an arbitrary $k^{\prime}$-subset of $V(P_{n})$. Then there must be a vertex
$v\in \{v_0,\ldots,v_{(n-1)/2-\lceil(k^{\prime}-1)/2\rceil}\}\cup \{v_{(n-1)/2+\lfloor(k^{\prime}-1)/2\rfloor+1},\ldots,v_{n-1}\}$ in $S^{\prime \prime}$. Thus,  $\varepsilon_{k,k^{\prime}}(S^{\prime \prime};P_n)\geq\max\{d_{P_n}(S^{\prime \prime}\cup \{v_0\}), d_{P_n}(S^{\prime \prime}\cup \{v_{n-1}\})\}\geq \lceil\frac{k^{\prime}-1}{2}\rceil+\frac{n-1}{2}$,
and we have $Sr_{k,k^{\prime}}(P_n)\geq \lceil\frac{k^{\prime}-1}{2}\rceil+\frac{n-1}{2}$. Hence, $Sr_{k,k^{\prime}}(P_n)= \lceil\frac{k^{\prime}-1}{2}\rceil+\frac{n-1}{2}$.
If $n$ is even, then there are two central vertices
$v_{\lfloor (n-1)/2\rfloor}$, $v_{\lceil (n-1)/2\rceil}$  in $P_{n}$. Let
$S^{\prime}=\{v_{\lfloor (n-1)/2\rfloor-\lceil(k^{\prime}-2)/2\rceil}, \ldots, v_{\lfloor (n-1)/2\rfloor}, v_{\lceil (n-1)/2\rceil}, \ldots, v_{\lceil (n-1)/2\rceil+\lfloor(k^{\prime}-2)/2\rfloor}\}$ and $S=S^{\prime}\cup \{v_{n-1}\}$.
Then
$Sr_{k,k^{\prime}}(P_n)\leq\varepsilon_{k,k^{\prime}}(S^{\prime};P_n)=d_{P_n}(S)=\lceil\frac{k^{\prime}-2}{2}\rceil+\frac{n}{2}$. Similarly, $Sr_{k,k^{\prime}}(P_n)\geq\lceil\frac{k^{\prime}-2}{2}\rceil+\frac{n}{2}$.
Hence, $Sr_{k,k^{\prime}}(P_n)=\lceil\frac{k^{\prime}-2}{2}\rceil+\frac{n}{2}$.

If $k-k^{\prime}\geq2$, then for any $k^{\prime}$-subset $S^{\prime}$, there exists a $k$-subset $S$ of $V(P_n)$ such that $S^{\prime}\subseteq S$ and $v_0$, $v_{n - 1}\in S$, then $\varepsilon_{k,k^{\prime}}(S^{\prime};P_n)\geq d_{P_n}(S)=n-1$. Moreover, since $\varepsilon_{k,k^{\prime}}(S^{\prime};P_n)\leq n-1$, we have $\varepsilon_{k,k^{\prime}}(S^{\prime};P_n)= n-1$. Hence, $Sr_{k,k^{\prime}}(P_n)=n-1$.

(3) Let $G=K_{n_1,n_2,\ldots,n_r}$ and $V(G)=V_1\cup V_2\cup \cdots \cup V_r$, where $|V_i|=n_i $ $(1\leq i\leq r)$.

If $k^{\prime}=1$, then $Sr_{k,1}(G)=\min\{\varepsilon_{k,1}(v;G): v\in V(G)\}$. For any $v\in V_i$, we have $\varepsilon_{k,1}(v;G)=k$ if and only if $|V_i|\geq k$; $\varepsilon_{k,1}(v;G)=k-1$ if and only if $|V_i|\leq k-1$, where $1\leq i\leq r$. Thus $Sr_{k,1}(G)=k$ if $n_1\geq k$, and $Sr_{k,1}(G)=k-1$ if $n_1\leq k-1$.

If $k^{\prime}\geq2$, let $S^{\prime}$ be a $k^{\prime}$-subset of $V(G)$ such that $S^{\prime}\cap V_1\neq \emptyset$ and $S^{\prime}\cap V_2\neq \emptyset$. Then, we have
$\varepsilon_{k,k^{\prime}}(S^{\prime};G)=k-1$, and hence $Sr_{k,k^{\prime}}(G)\leq k-1$. Moreover, since $Sr_{k,k^{\prime}}(G)\geq k-1$, we have $Sr_{k,k^{\prime}}(G)=k-1$. $\Box$

The following theorem derives a sharp upper bound of the Steiner $k$-diameter by the Steiner $k$-radius of a tree.

\begin{theorem}(\cite{Steiner distance})\label{lem2.6}
Let $n$, $k$ be two integers with $n\geq k\geq2$, and $T$ be a tree with order $n$. Then
$$Sd_k(T)\leq \frac{k}{k-1}Sr_{k}(T),$$
the equality holds if $T\cong K_{1,n-1}$.
\end{theorem}

Henning, Oellermann and Swart \cite{diameter radius1} put forward Conjecture \ref{con2.7} and proved that it holds for $k=3$ and $k=4$.

\begin{conjecture}(\cite{diameter radius1})\label{con2.7}
Let $k\geq 2$ be an integer. If $G$ is a connected graph with order at least $k$. Then
\begin{align*}
Sd_k(G)\leq\dfrac{2(k+1)}{2k-1}Sr_{k}(G).
\end{align*}
\end{conjecture}


Subsequently, Reiswig \cite{diameter radius2} improved this bound in Conjecture \ref{con2.7} for $k\geq 5$ and obtained the following tight bound.

\begin{theorem}(\cite{diameter radius2})\label{lem2.8}
Let $k\geq 5$ be an integer. If $G$ is a connected graph with order at least $k$. Then
\begin{align*}
Sd_k(G)\leq \dfrac{k+3}{k+1}Sr_{k}(G).
\end{align*}
Moreover, the upper bound is tight.
\end{theorem}

\section{\boldmath Steiner $k$-diameter and Steiner $(k,k^{\prime})$-radius of trees}

~~~~Let $T$ be a tree with $|L(T)|\geq 3$, and let $S$ be a set of pendant vertices with $|S|\geq 3$. For a vertex $v\in S$, $l_{v}(S)$ denotes the shortest distance from $v$ to a branching vertex in $T_S$. The following lemma plays an important role in the proof of the main results in this section.
\begin{lemma}\label{lem3.1}(\cite{Steiner distance})
Let $k\geq 3$ be an integer, and $T$ be a tree with at least $k$ pendant vertices.
Let $S$ be a $k$-subset of $L(T)$ and $v\in S$. Then
$d_T(S)=d_T(S\setminus \{v\})+l_{v}(S)$.
\end{lemma}

Next, we will give an upper bound of the Steiner $k$-diameter of a tree by its Steiner $(k-k^{\prime})$-diameter.

\begin{theorem}\label{thm3.2}
Let $n$, $k$, $k^{\prime}$ be integers with $n\geq k\geq 3$, $k^{\prime}\geq 1$ and $k-k^{\prime}\geq 2$. Let $T$ be a tree of order $n$. Then
$$Sd_{k}(T)\leq \frac{k}{k-k^{\prime}} Sd_{k-k^{\prime}}(T),$$
where the equality holds if $T$ is a starlike tree with at least $k$ longest pendant paths.
\end{theorem}

\noindent{\bf Proof.} Let $S\subseteq V(T)$ such that $d_{T}(S)=Sd_{k}(T)$, where $|S|=k$. We will consider two cases based on the number of pendant vertices of the tree $T$.

\noindent{\bf Case 1. } $|L(T)|\leq k-k^{\prime}$.

If $|L(T)|\leq k-k^{\prime}$, then $Sd_{k}(T)=Sd_{k-k^{\prime}}(T)=n-1$.

\noindent{\bf Case 2. } $ |L(T)|\geq k-k^{\prime}+1$.

Let $S_1=S\cap L(T)$ and $s=|S_1|$. Then $S_1\subseteq L(T)$ and $Sd_{k}(T)=d_{T}(S)=d_{T}(S_1)$, where $k-k^{\prime}+1\leq s\leq k$.

In the following, we will show that there is at least one vertex $v\in S_1$, such that $l_{v}(S_1)\leq \frac{1}{k-k^{\prime}}Sd_{k-k^{\prime}}(T)$.
Assume, by contrary, that for any vertex $u\in S_1$, it always holds that $l_{u}(S_1)> \frac{1}{k-k^{\prime}}Sd_{k-k^{\prime}}(T)$. Let $S^{\prime}$ be the set of $k-k^{\prime}$ vertices of $S_1$. Then
\begin{align*}
d_{T}(S^{\prime})\geq \sum_{u\in S^{\prime}}l_u(S_1)>(k-k^{\prime})\frac{1}{k-k^{\prime}}Sd_{k-k^{\prime}}(T)=Sd_{k-k^{\prime}}(T).
\end{align*}
This result contradicts the definition of the Steiner $(k-k^{\prime})$-diameter. Hence, there must be at least one vertex $v\in S_1$ such that $l_{v}(S_1)\leq \frac{1}{k-k^{\prime}}Sd_{k-k^{\prime}}(T)$.
Let $S_2=S_1\setminus\{v_1\}$, where $v_1\in S_1$ and $l_{v_1}(S_1)\leq \frac{1}{k-k^{\prime}}Sd_{k-k^{\prime}}(T)$. By Lemma \ref{lem3.1}, we have
\begin{align*}
d_{T}(S_1)= d_{T}(S_2)+l_{v_{1}}(S_1)\leq d_{T}(S_2)+\frac{1}{k-k^{\prime}}Sd_{k-k^{\prime}}(T).
\end{align*}

Iteratively, we can derive that the set $S_{s- k + k^{\prime}}$ satisfies $|S_{s-k+ k^{\prime}}|=k-k^{\prime}+1\geq 3$, and there exists a vertex $v_{s-k+k^{\prime}}$ in $S_{s-k+k^{\prime}}$ such that $l_{v_{s-k+k^{\prime}}}(S_{s-k+k^{\prime}})\leq \frac{1}{k-k^{\prime}}Sd_{k-k^{\prime}}(T)$. By Lemma \ref{lem3.1}, we have
\begin{align*}
d_T(S_{s-k+k^{\prime}})=d_T(S_{s-k+k^{\prime}+1})+l_{v_{s-k+k^{\prime}}}(S_{s-k+k^{\prime}})\leq d_T(S_{s-k+k^{\prime}+1})+\frac{1}{k-k^{\prime}}Sd_{k-k^{\prime}}(T),
\end{align*}
where $S_{s-k+k^{\prime}+1}=S_{s-k+k^{\prime}}\setminus \{v_{s-k+k^{\prime}}\}$ and $|S_{s-k+k^{\prime}+1}|=k-k^{\prime}$.
Thus, we have
\begin{align*}
d_{T}(S_1)&\leq d_{T}(S_2)+\frac{1}{k-k^{\prime}}Sd_{k-k^{\prime}}(T)\\
&\leq d_{T}(S_{s-k+k^{\prime}+1})+\frac{s-k+k^{\prime}}{k-k^{\prime}}Sd_{k-k^{\prime}}(T)\\
&=\frac{s}{k-k^{\prime}}Sd_{k-k^{\prime}}(T)\\
&\leq \frac{k}{k-k^{\prime}}Sd_{k-k^{\prime}}(T).
\end{align*}

Therefore, we conclude that
$Sd_{k}(T)\leq \frac{k}{k-k^{\prime}}Sd_{k-k^{\prime}}(T)$.

If $T$ is a starlike tree with at least $k$ pendant paths with the longest length $l$, then $Sd_k(T)=kl$ and $Sd_{k-k^{\prime}}(T)=(k-k^{\prime})l$. Hence, $Sd_k(T)=\dfrac{k}{k-k^{\prime}}Sd_{k-k^{\prime}}(T)$. $\Box$

\begin{theorem}\label{thm3.3}
Let $n$, $k$, $k^{\prime}$ be integers with $n\geq k\geq 3$, $k^{\prime}\geq 1$ and $k-k^{\prime}\geq 2$. Let $T$ be a tree of order $n$. Then
$$Sd_{k-k^{\prime}}(T)\leq Sr_{k,k^{\prime}}(T).$$
\end{theorem}
\noindent{\bf Proof.} Suppose $Sr_{k,k^{\prime}}(T)=\varepsilon_{k,k^{\prime}}(S^{\prime};T)=d_{T}(S)$, where $S^{\prime}\subseteq S\subseteq V(T)$, $|S^{\prime}|=k^{\prime}$ and $|S|=k$.
It is sufficient to show that $d_{T}(S)\geq d_{T}(S^{\prime \prime})$ for any $(k-k^{\prime})$-subset $S^{\prime \prime}$ in $T$.

Assume that there exists a $(k-k^{\prime})$-subset $S^{\prime \prime}$ of $V(T)$ such that $d_{T}(S)< d_{T}(S^{\prime \prime})$. It follows that the inequality $d_{T}(S)<d_{T}(S^{\prime \prime})\leq d_{T}(S^{\prime}\cup S^{\prime \prime})$ holds.

If $|S^{\prime}\cup S^{\prime \prime}|=k$, then we get another $k$-subset $S^{\prime}\cup S^{\prime \prime}$ satisfying $S^{\prime}\subseteq S^{\prime}\cup S^{\prime \prime}$ and $d_{T}(S)<d_{T}(S^{\prime}\cup S^{\prime \prime})$, a contradiction.

If $|S^{\prime}\cup S^{\prime \prime}|=k^{\prime \prime}<k$, then by taking any $k-k^{\prime \prime}$ vertices of $V(T)\setminus(S^{\prime}\cup S^{\prime \prime})$ and adding them to $S^{\prime}\cup S^{\prime \prime}$, we can construct a new $k$-subset $S^{*}$ satisfying $S^{\prime}\subseteq S^{*}$ and
$$d_{T}(S)<d_{T}(S^{\prime}\cup S^{\prime \prime})\leq d_{T}(S^*),$$
which is also a contradiction.

Consequently, we have $d_{T}(S)\geq d_{T}(S^{\prime \prime})$ for any $S^{\prime \prime}\subseteq V(T)$ with $|S^{\prime \prime}|=k-k^{\prime}$. Therefore, we have $Sr_{k,k^{\prime}}(T)=d_{T}(S)\geq Sd_{k-k^{\prime}}(T)$. $\Box$

\begin{theorem}\label{thm3.4}
Let $n$, $k$, $k^{\prime}$ be integers with $n\geq k\geq 3$ and $k>k^{\prime}\geq 1$. Let $T$ be a tree with order $n$. Then
$$Sd_{k}(T)\leq \frac{k}{k-k^{\prime}}Sr_{k,k^{\prime}}(T).$$
\end{theorem}
\noindent{\bf Proof.} Denote $d=diam(T)$. By the proof of Corollary \ref{cor2}, we have $Sd_{k}(T)=d(S)\leq d+(k-2)\lfloor\frac{d}{2}\rfloor\leq k\frac{d}{2}$, where $T_S$ contains a diametrical path $P_{d+1}=u_0u_1\cdots u_d$ of $T$.
Since $Sr_{2,1}(T) \geq \frac{d}{2}$, it follows that $Sd_k(T)\leq kSr_{2,1}(T)$.

If $k-k^{\prime}=1$, then for any $(k-1)$-subset $S^{\prime}\subseteq V(T)$, we have $\varepsilon_{k,k-1}(S^{\prime};T)\geq \max\{d_T(S^{\prime}\cup \{u_0\}),d_T(S^{\prime}\cup \{u_d\})\}\geq \varepsilon_{2,1}(v;T)\geq Sr_{2,1}(T)$, where $v\in S^{\prime}$. Therefore, we have $Sr_{k,k-1}(T)\geq Sr_{2,1}(T)$ and $Sd_{k}(T)\leq kSr_{k,k-1}(T)$.
If $k-k^{\prime}\geq2$, then Theorems \ref{thm3.2} and \ref{thm3.3} imply that $Sd_{k}(T)\leq \frac{k}{k-k^{\prime}}Sr_{k,k^{\prime}}(T).$

Hence, for $k\geq 3$, $k>k^{\prime}\geq 1$, we have $Sd_{k}(T)\leq \frac{k}{k-k^{\prime}}Sr_{k,k^{\prime}}(T).$ $\Box$

\section{\boldmath Steiner $k$-diameter and Steiner $(k,k^{\prime})$-radius of trees for $k^{\prime}=2,3$}
~~~~In the preceding section, we obtained an upper bound of the Steiner $k$-diameter of a tree by its Steiner ($k,k^{\prime}$)-radius for $k\geq3$ and $k>k^{\prime}\geq1$. Nevertheless, for $k^{\prime}>1$, it must be noted that the bound is not tight. Therefore, in this section, we will focus on the special cases of $k^{\prime}=2, 3$ and give a tight upper bound of the Steiner $k$-diameter of a tree by its Steiner $(k,2)$-radius and Steiner $(k,3)$-radius.

Let $P_{l}$ denote the path on $l$ vertices. The graph $P_{l}(a,b;x)$ is obtained from $P_l$ by attaching $a$ and $b$ pendant paths of length $x$ to the two end vertices of $P_l$, respectively, where $a$, $b$, $l$, $x$ are positive integers.

\begin{definition}\label{def4.1}
Let $T$ be a tree and $P=u_0 u_1\cdots u_d$ be a diametrical path of $T$. Define $A=A(T,P)$ as the maximal set of pendant vertices of $T$ satisfying the following three conditions:
\begin{enumerate}[label=(\arabic*)]
  \item $u_0, u_d\in A$;
  \item the distance of each vertex in $A$ to $C(T)$ is $\left\lfloor \frac{d}{2} \right\rfloor$;
  \item the paths $\{ P_w: w \in A \}$ are mutually edge-disjoint, where $P_w$ is the shortest path from $w$ to $C(T)$ for each $w\in A$.
\end{enumerate}
\end{definition}

\begin{lemma}\label{lem4.0}
Let $T$ be a tree and $P=u_0 u_1\cdots u_d$ be a diametrical path of $T$.
Let $A=A(T,P)$ be the subset of $L(T)$ defined in Definition \ref{def4.1}.
Then $2\leq |A|\leq |L(T)|$ and $d_{T}(A)=Sd_{|A|}(T)$.
\end{lemma}
\noindent{\bf Proof.} Since $u_0, u_d\in A$ and $A\subseteq L(T)$, we have $2\leq |A|\leq |L(T)|$. By the definition of $A$, we have $d_T(A) = d+(|A|-2)\lfloor\frac{d}{2}\rfloor$. By Corollary \ref{cor2}, $Sd_{|A|}(T)\leq d+(|A|-2)\lfloor\frac{d}{2}\rfloor$. Hence, $Sd_{|A|}(T)\leq d_T(A)$.
On the other hand, $d_T(A)\leq Sd_{|A|}(T)$. Therefore, we obtain that $ d_T(A)=Sd_{|A|}(T)$. $\Box$

\begin{lemma}\label{lem4.1}
Let $n$, $k$ be two integers with $n\geq k\geq 3$ and $T$ be a tree of order $n$. Let $P=u_0 u_1\cdots u_d$ be a diametrical path of $T$ and $S^{\prime}=\{u_{\lfloor d/2\rfloor},u_{\lfloor d/2\rfloor+1}\}$. Then
$$\varepsilon_{k,2}(S^{\prime};T)=Sr_{k,2}(T).$$
\end{lemma}
\noindent{\bf Proof.}
Let $A=A(T,P)$ be the subset of $L(T)$ defined in Definition \ref{def4.1}. We consider two cases according to the value of $|A|$.

\noindent{\bf Case 1. } $|A|\leq k-2$.

By $2\leq |A|\leq k-2$, we have $k\geq 4$. By Theorem \ref{thm3.3}, $Sr_{k,2}(T)\geq Sd_{k-2}(T)$.
If $|L(T)|\leq k-2$, then $\varepsilon_{k,2}(S^{\prime};T)=|E(T)|=Sd_{k-2}(T)$. If $|L(T)|>k-2$, then $|L(T)\setminus S^{\prime}|=|L(T)|>k-2>0$.
By Lemma \ref{steinertree}, given a Steiner $(k-1)$-ecc $S^{\prime}$-tree $T^{k-1}$, there exists a Steiner $k$-ecc $S^{\prime}$-tree $T^{k}$ such that $T^{k}$ contains $T^{k-1}$. By the definition of the Steiner $(k,k^{\prime})$-eccentricity, we have
\[
\varepsilon_{|A|+2,2}(S^{\prime};T)=\max\{d_{T}(S^{\prime}\cup S): S\subseteq V(T)\setminus S^{\prime}, |S|=|A|\}=d_{T}(S^{\prime}\cup A).
\]
Thus, for the Steiner $(|A|+2)$-ecc $S^{\prime}$-tree $T_{S^{\prime}\cup A}$, there exists a Steiner $k$-ecc $S^{\prime}$-tree $T_{S^{\prime}\cup S^*}$ such that $T_{S^{\prime}\cup S^*}$ contains $T_{S^{\prime}\cup A}$, where $|S^*| = k - 2$. It follows that $A\subseteq S^*$, and hence $T_A \subseteq T_{S^*}$. Since $u_0, u_d \in A \subseteq S^*$, we have $S^{\prime}\subseteq V(T_{S^*})$. Therefore,
\[
d_{T}(S^{\prime} \cup S^*) = d_{T}(S^*),
\]
which implies
\[
\varepsilon_{k, 2}(S^{\prime}; T) = d_{T}(S^*) \leq Sd_{k-2}(T).
\]
By Theorem \ref{thm3.3}, $\varepsilon_{k,2}(S^{\prime};T)\geq Sr_{k,2}(T)\geq Sd_{k-2}(T)$. It follows that $\varepsilon_{k,2}(S^{\prime};T)=Sd_{k-2}(T)$ and $\varepsilon_{k,2}(S^{\prime};T)=Sr_{k,2}(T)$.

\noindent{\bf Case 2 } $|A|\geq k-1$.

By $|A|\geq k-1$, we have $\varepsilon_{k,2}(S^{\prime};T)=d_{T}(S^{\prime}\cup A^{*})=(k-2)\lfloor\frac{d}{2}\rfloor+1$, where $A^{*}\subseteq A\setminus\{u_d\}$ and $|A^{*}|=k-2$. In the following, we will show that $\varepsilon_{k,2}(S;T)\geq (k-2)\lfloor\frac{d}{2}\rfloor+1$ for any $S=\{a,b\}\subseteq V(T)$.

If $a$, $b\in V(T_A)$, then $\varepsilon_{k,2}(S;T)\geq \varepsilon_{k,2}(S;T_A)= \max\{d_{T_A}(S\cup A^{\prime}): A^{\prime}\subseteq A, |A^{\prime}|=k-2\}\geq (k-2)\lfloor\frac{d}{2}\rfloor+1$.
Suppose at least one of $a$ and $b$ is not in $V(T_A)$. Without loss of generality, assume $a\notin V(T_A)$. Then we have
$\varepsilon_{k,2}(S;T)\geq d_{T}(S\cup A^{\prime \prime})\geq d_{T}(A^{\prime \prime})+d_T(a,T_{A^{\prime \prime}})\geq (k-2)\lfloor\frac{d}{2}\rfloor+1$, where $A^{\prime \prime}\subseteq A$ and $|A^{\prime \prime}|=k-2$.

Thus, we obtain that $\varepsilon_{k,2}(S;T)\geq \varepsilon_{k,2}(S^{\prime};T)$ for any $S=\{a,b\}\subseteq V(T)$. It follows that $\varepsilon_{k,2}(S^{\prime};T)=Sr_{k,2}(T)$.

The proof is thus complete. $\Box$

\begin{theorem}\label{them4.3}
Let $n$, $k$ be two integers with $n\geq k\geq 3$ and $T$ be a tree of order $n$. Then
$$Sd_{k}(T)\leq \frac{k}{k-2}Sr_{k,2}(T)-\frac{2}{k-2},$$
the equality holds if $T\cong P_{2}(a,b;x)$, where $\frac{n-2}{x}$ is an integer and $a+b=\frac{n-2}{x}\geq k$.
\end{theorem}

\noindent{\bf Proof.} Let $P=u_0 u_1\cdots u_d$ be a diametrical path of $T$, and $A=A(T,P)$ be the subset of $L(T)$ defined in Definition \ref{def4.1}. Let $S^{\prime}=\{u_{\lfloor d/2\rfloor},u_{\lfloor d/2\rfloor+1}\}$. By Lemma \ref{lem4.1},
$Sr_{k,2}(T)=\varepsilon_{k,2}(S^{\prime};T)$.

If $|L(T)|\leq k-1$, then $Sd_{k}(T)=Sd_{k-1}(T)=|E(T)|$. By Theorem \ref{thm3.2}, $Sd_{k-1}(T)\leq \frac{k-1}{k-2}Sd_{k-2}(T)$. By Theorem \ref{thm3.3}, $Sd_{k-2}(T)\leq Sr_{k,2}(T)$.
Hence, we obtain that $Sd_{k-1}(T)\leq \frac{k-1}{k-2}Sr_{k,2}(T)$. By $k\geq3$, we have $Sr_{k,2}(T)\geq k-1\geq 2$. Therefore,
$$Sd_{k}(T)=Sd_{k-1}(T)\leq\frac{k-1}{k-2}Sr_{k,2}(T)\leq \frac{k}{k-2}Sr_{k,2}(T)- \frac{2}{k-2}.$$

Assume $|L(T)| \geq k$ in the following proof. We consider two cases based on the value of $|A|$.

\noindent{\bf Case 1. } $|A|\leq k-2$.

By Lemma \ref{lem4.0}, $d_T(A)=Sd_{|A|}(T)$.
Since $|L(T)|\geq k$ and $|A|\leq k-2$,
by Lemma \ref{steinertree}, for $p = k-2, k-1, k $, there exist $S_p$ such that $d_T(S_p)=Sd_p(T)$ and $T_{A}\subseteq T_{S_{k-2}} \subset T_{S_{k-1}} \subset T_{S_k}$.
Thus, we obtain that $A\subseteq S_{k-2}\subset S_{k-1}\subset S_{k}$.
Without loss of generality, let $S_{k-2} = \{v_1, v_2, \ldots, v_{k-2}\}$, $S_{k-1} = \{v_1, v_2, \ldots, v_{k-1}\}$ and $S_k = \{v_1, v_2, \ldots, v_k\}$.
Then we have $$Sd_k(T)=Sd_{k-2}(T)+l_{v_{k-1}}(S_{k-1})+l_{v_{k}}(S_k).$$
We now prove $l_{v_k}(S_k)\leq l_{v_{k-1}}(S_{k-1})$. Suppose, to the contrary, that $l_{v_k}(S_k)> l_{v_{k-1}}(S_{k-1})$. Then $l_{v_k}(S_k)> l_{v_{k-1}}(S_{k-1})\geq l_{v_{k-1}}(S_k)$. By Lemma \ref{lem3.1}, $d_{T}(S_{k-1})=d_{T}(S_k)-l_{v_{k}}(S_k)$. It follows that
$$d_{T}(S_{k-1})< d_{T}(S_k)-l_{v_{k-1}}(S_k)=d_{T}(S_k\setminus\{v_{k-1}\}).$$
This contradicts the choice of $S_{k-1}$. Thus, $l_{v_k}(S_k)\leq l_{v_{k-1}}(S_{k-1})$.
Similarly, it follows that $l_{v_{k-1}}(S_{k-1})\leq l_{v}(S_{k-2})$ for any $v\in S_{k-2}$.
Thus, we have $Sd_k(T)\leq Sd_{k-2}(T)+2l_{v_{k-1}}(S_{k-1})$ and
\begin{align*}
Sd_{k-2}(T)=d_{T}(S_{k-2})&\geq d_T(A)+\sum_{v\in S_{k-2}\setminus A}l_v(S_{k-2})\\
&\geq d+(|A|-2)\lfloor\frac{d}{2}\rfloor+(k-2-|A|)l_{v_{k-1}}(S_{k-1}).
\end{align*}
Since $A\subseteq S_{k-2}$ and $v_{k-1}\notin S_{k-2}$, we have $l_{v_{k-1}}(S_{k-1})\leq \lfloor\frac{d}{2}\rfloor-1$.
Then
\begin{align*}
Sd_{k}(T)-\frac{k}{k-2}Sr_{k,2}(T)+\frac{2}{k-2}&\leq Sd_{k-2}(T)+2l_{v_{k-1}}(S_{k-1})-\frac{k}{k-2}Sd_{k-2}(T)+\frac{2}{k-2}\\
&=
-\frac{2}{k-2}Sd_{k-2}(T)+2l_{v_{k-1}}(S_{k-1})+\frac{2}{k-2}\\
&\leq -\frac{2}{k-2}(d+(|A|-2)\lfloor\frac{d}{2}\rfloor+(k-2-|A|)l_{v_{k-1}}(S_{k-1}))\\
&\quad +2l_{v_{k-1}}(S_{k-1})+\frac{2}{k-2}\\
&= -\frac{2}{k-2}(d+(|A|-2)\lfloor\frac{d}{2}\rfloor)+\frac{2|A|}{k-2}l_{v_{k-1}}(S_{k-1})+\frac{2}{k-2}\\
&\leq -\frac{2}{k-2}(d+(|A|-2)\lfloor\frac{d}{2}\rfloor)+\frac{2|A|}{k-2}(\lfloor\frac{d}{2}\rfloor-1)+\frac{2}{k-2}\\
&\leq -\frac{2}{k-2}(d-2\lfloor\frac{d}{2}\rfloor)-\frac{2|A|}{k-2}+\frac{2}{k-2}\\
&\leq -\frac{2|A|}{k-2}+\frac{2}{k-2}<0.
\end{align*}
Thus, we obtain that $Sd_{k}(T)<\frac{k}{k-2}Sr_{k,2}(T)-\frac{2}{k-2}$.

\noindent{\bf Case 2. } $|A|\geq k-1$.

By Corollary \ref{cor2}, $Sd_k(T)\leq d+(k-2)\lfloor\frac{d}{2}\rfloor$.
By $|A|\geq k-1$, we have $Sr_{k,2}(T)=\varepsilon_{k,2}(S^{\prime};T)=(k-2)\lfloor\frac{d}{2}\rfloor+1$. It follows that
\begin{align*}
Sd_{k}(T)-\frac{k}{k-2}Sr_{k,2}(T)+\frac{2}{k-2}\leq d-2\lfloor\frac{d}{2}\rfloor-1\leq0.
\end{align*}
Thus, we obtain that $Sd_{k}(T)\leq\frac{k}{k-2}Sr_{k,2}(T)-\frac{2}{k-2}$.

If $T\cong P_{2}(a,b;x)$, then $Sd_{k}(T)=kx+1$ and $Sr_{k,2}(T)=(k-2)x+1$. Thus, $Sd_k(T)=\frac{k}{k-2}Sr_{k,2}(T)-\frac{2}{k-2}$.

Therefore, we conclude that $Sd_k(T)\leq\frac{k}{k-2}Sr_{k,2}(T)-\frac{2}{k-2}$, where the equality holds if $T\cong P_{2}(a,b;x)$. $\Box$

\begin{lemma}\label{lem4.2}
Let $n$, $k$ be two integers with $n\geq k\geq 4$ and $T$ be a tree of order $n$. Let $P=u_0 u_1\cdots u_d$ be a diametrical path of $T$ and $S^{\prime \prime}=\{u_{\lfloor d/2\rfloor-1},u_{\lfloor d/2\rfloor},u_{\lfloor d/2\rfloor+1}\}$. Then
$$\varepsilon_{k,3}(S^{\prime\prime};T)=Sr_{k,3}(T).$$
\end{lemma}
\noindent{\bf Proof.}
Let $A=A(T,P)$ be the subset of $L(T)$ defined in Definition \ref{def4.1}.
We consider two cases according to the value of $|A|$.

\noindent{\bf Case 1. } $|A|\leq k-3$.

By $2\leq |A|\leq k-3$, we have $k\geq 5$. By Theorem \ref{thm3.3}, $Sr_{k,3}(T)\geq Sd_{k-3}(T)$.
If $|L(T)|\leq k-3$, then $\varepsilon_{k,3}(S^{\prime \prime};T)=|E(T)|=Sd_{k-3}(T)$. If $|L(T)|>k-3$, then $|L(T)\setminus S^{\prime \prime}|=|L(T)|>k-3>0$. By the definition of the Steiner $(k,k^{\prime})$-eccentricity, we have $\varepsilon_{|A|+3,3}(S^{\prime \prime};T)=d_{T}(S^{\prime \prime}\cup A)$. By Lemma \ref{steinertree}, for the tree $T_{S^{\prime \prime}\cup A}$, there exists a Steiner $k$-ecc $S^{\prime \prime}$-tree $T_{S^{\prime \prime}\cup S^*}$ such that $T_{S^{\prime \prime}\cup A}\subseteq T_{S^{\prime \prime}\cup S^*}$, where $|S^*|=k-3$. It follows that $A \subseteq S^*$.
Since $u_0, u_d \in A \subseteq S^*$, we have $S^{\prime \prime}\subseteq V(T_{S^*})$. Therefore,
\[
\varepsilon_{k, 3}(S^{\prime \prime}; T) =d_{T}(S^{\prime \prime} \cup S^*) = d_{T}(S^*) \leq Sd_{k-3}(T).
\]
By Theorem \ref{thm3.3}, $\varepsilon_{k,3}(S^{\prime\prime};T)\geq Sr_{k,3}(T)\geq Sd_{k-3}(T)$. It follows that $\varepsilon_{k,3}(S^{\prime\prime};T)= Sd_{k-3}(T)$ and $\varepsilon_{k,3}(S^{\prime\prime};T)= Sr_{k,3}(T)$.

\noindent{\bf Case 2. } $|A|\geq k-2$.

If $|A|=k-2$, then $\varepsilon_{k,3}(S^{\prime\prime};T)=d_{T}(S^{\prime\prime}\cup A\setminus\{u_0\})=(k-4)\lfloor\frac{d}{2}\rfloor+\lceil\frac{d}{2}\rceil+1$. If $|A|\geq k-1$, then $\varepsilon_{k,3}(S^{\prime\prime};T)=d_{T}(S^{\prime\prime}\cup A^*)=(k-3)\lfloor\frac{d}{2}\rfloor+2$, where $A^*\subseteq A\setminus\{u_0,u_d\}$ and $|A^*|=k-3$. In the following, we will show that $\varepsilon_{k,3}(S;T)\geq \varepsilon_{k,3}(S^{\prime \prime};T)$ for any $S=\{a,b,c\}\subseteq V(T)$.

Suppose $a,b,c\in V(T_A)$.
If $|A|=k-2$, then $\varepsilon_{k,3}(S;T)\geq \varepsilon_{k,3}(S;T_A)\geq d_{T}(A\setminus \{u_0\})+1=(k-4)\lfloor\frac{d}{2}\rfloor+\lceil\frac{d}{2}\rceil+1$. If $|A|\geq k-1$, then $\varepsilon_{k,3}(S;T)\geq\varepsilon_{k,3}(S;T_A)\geq \sum_{w\in A_1}|E(P_w)|+2=(k-3)\lfloor\frac{d}{2}\rfloor+2$, where $A_1\subseteq A$ and $|A_1|=k-3$.
%
%
Suppose exactly one of the vertices $a$, $b$ and $c$ is not in $V(T_A)$, without loss of generality, assume $a\notin V(T_A)$. Then
$\varepsilon_{k,3}(S;T)\geq \varepsilon_{k-1,2}(\{b,c\};T_A)+1\geq (k-3)\lfloor\frac{d}{2}\rfloor+2$. Suppose at least two of the vertices $a$, $b$ and $c$ are not in $V(T_A)$, without loss of generality, assume $a,b \notin V(T_A)$. Then $\varepsilon_{k,3}(S;T)\geq d_{T}(\{a,b\} \cup A_2)\geq d_{T}(A_2)+2 \geq (k-3)\lfloor\frac{d}{2}\rfloor+2$, where $A_2\subseteq A$ and $|A_2|=k-3$.
Thus, for any vertex set $S = \{a, b, c\} \subseteq V(T)$, if $|A| = k-2 $, then $\varepsilon_{k,3}(S;T)\geq (k-4)\lfloor\frac{d}{2}\rfloor+\lceil\frac{d}{2}\rceil+1$; if $|A| \geq k-1$, then $\varepsilon_{k,3}(S;T)\geq (k-3)\lfloor\frac{d}{2}\rfloor+2$.

Thus, we obtain that $\varepsilon_{k,3}(S;T)\geq \varepsilon_{k,3}(S^{\prime\prime};T)$ for any $S=\{a,b,c\}\subseteq V(T)$. It follows that $\varepsilon_{k,3}(S^{\prime\prime};T)=Sr_{k,3}(T)$. $\Box$

\begin{theorem}\label{thm4.4}
Let $n$, $k$ be two integers with $n\geq k\geq 4$ and $T$ be a tree of order $n$. Then
$$Sd_{k}(T)\leq \frac{k}{k-3}Sr_{k,3}(T)-\frac{6}{k-3},$$
the equality holds if $T\cong P_{3}(a,b;x)$, where $\frac{n-3}{x}$ is an integer and $a+b= \frac{n-3}{x}\geq k$.
\end{theorem}

\noindent{\bf Proof.}  Let $P=u_0 u_1\cdots u_d$ be a diametrical path of $T$, and $A=A(T,P)$ be the subset of $L(T)$ defined in Definition \ref{def4.1}. Let $S^{\prime \prime}=\{u_{\lfloor d/2\rfloor-1},u_{\lfloor d/2\rfloor},u_{\lfloor d/2\rfloor+1}\}$. By Lemma \ref{lem4.2}, we have
$Sr_{k,3}(T)=\varepsilon_{k,3}(S^{\prime \prime};T)$.

If $|L(T)|\leq k-2$, then $Sd_{k}(T)=Sd_{k-2}(T)=|E(T)|$. By Theorem \ref{thm3.2}, $Sd_{k-2}(T)\leq \frac{k-2}{k-3}Sd_{k-3}(T)$. By Theorem \ref{thm3.3}, $Sd_{k-3}(T)\leq Sr_{k,3}(T)$. Hence, we obtain that $Sd_{k-2}(T)\leq \frac{k-2}{k-3}Sr_{k,3}(T)$. By $k\geq 4$, we have $Sr_{k,3}(T)\geq k-1\geq 3$. Therefore,
$$Sd_{k}(T)\leq\frac{k-2}{k-3}Sr_{k,3}(T) \leq \frac{k}{k-3}Sr_{k,3}(T)-\frac{6}{k-3}.$$

Now it suffices to prove that the result holds for $|L(T)|\geq k-1$. We consider three cases based on the value of $|A|$.

\noindent{\bf Case 1. } $|A|\leq k-3$.

If $|L(T)|=k-1$, then $Sd_{k}(T)=Sd_{k-1}(T)=|E(T)|$. By Lemma \ref{lem4.0}, $d_{T}(A)=Sd_{|A|}(T)$.
By Lemma \ref{steinertree}, $T_A\subseteq T_{S_{k-3}}\subset T_{S_{k-2}}\subset T_{S_{k-1}}$, where $S_p$ is a $p$-subset of $V(T)$ such that $d_T(S_p) = Sd_p(T)$ for $p=k-3, k-2, k-1$.
Thus, we obtain that $A\subseteq S_{k-3}\subset S_{k-2}\subset S_{k-1}$.
Without loss of generality, let $S_{k-3}=\{v_1,v_2,\ldots,v_{k-3}\}$, $S_{k-2}=\{v_1,v_2,\ldots,v_{k-2}\}$ and $S_{k-1}=\{v_1,v_2,\ldots,v_{k-1}\}$.
By analogy with the proof of Theorem \ref{them4.3}, we have $l_{v_{k-1}}(S_{k-1})\leq l_{v_{k-2}}(S_{k-2})\leq \lfloor\frac{d}{2}\rfloor-1$ and $l_{v_{k-2}}(S_{k-2})\leq l_v(S_{k-3})$ for any $v\in S_{k-3}$.
Thus,
$$Sd_{k-1}(T)=Sd_{k-3}(T)+l_{v_{k-2}}(S_{k-2})+l_{v_{k-1}}(S_{k-1})\leq Sd_{k-3}(T)+2l_{v_{k-2}}(S_{k-2})$$
and
\begin{align*}
Sd_{k-3}(T)=d_{T}(S_{k-3})&\geq d_T(A)+\sum_{v\in S_{k-3}\setminus A}l_v(S_{k-3})\\
&\geq d+(|A|-2)\lfloor\frac{d}{2}\rfloor+(k-3-|A|)l_{v_{k-2}}(S_{k-2}).
\end{align*}

Therefore,
\begin{align*}
Sd_{k}(T)-\frac{k}{k-3}Sr_{k,3}(T)+\frac{6}{k-3}&\leq Sd_{k-3}(T)+2l_{v_{k-2}}(S_{k-2})-\frac{k}{k-3}Sr_{k,3}(T)+\frac{6}{k-3}\\
&= -\frac{3}{k-3}Sd_{k-3}(T)+2l_{v_{k-2}}(S_{k-2})+\frac{6}{k-3}\\
&\leq  -\frac{3}{k-3}(d+(|A|-2)\lfloor\frac{d}{2}\rfloor+(k-3-|A|)l_{v_{k-2}}(S_{k-2}))\\
&\quad +2l_{v_{k-2}}(S_{k-2}) +\frac{6}{k-3}\\
&=-\frac{3}{k-3}(d+(|A|-2)\lfloor\frac{d}{2}\rfloor)+(\frac{3|A|}{k-3}-1)l_{v_{k-2}}(S_{k-2}) +\frac{6}{k-3}\\
&\leq  -\frac{3}{k-3}(d+(|A|-2)\lfloor\frac{d}{2}\rfloor)+(\frac{3|A|}{k-3}-1)(\lfloor\frac{d}{2}\rfloor-1) +\frac{6}{k-3}\\
&\leq -\frac{3}{k-3}(d-2\lfloor\frac{d}{2}\rfloor)-\frac{3|A|}{k-3}+\frac{6}{k-3}-\lfloor\frac{d}{2}\rfloor+1\\
&\leq -\frac{3|A|}{k-3}+\frac{6}{k-3}\leq0.
\end{align*}
Thus, we obtain that $Sd_{k}(T)\leq\frac{k}{k-3}Sr_{k,3}(T)-\frac{6}{k-3}$.


If $|L(T)|\geq k$, then by Lemma \ref{lem2}, we obtain a $k$-subset $S_{k}$ such that $S_{k-1}\subset S_k$ and $d_T(S_k)= Sd_k(T)$. Let $S_k=\{v_1,v_2,\ldots,v_k\}$.
Then we have $l_{v_{k}}(S_k)\leq l_{v_{k-1}}(S_{k-1})\leq l_{v_{k-2}}(S_{k-2})\leq \lfloor\frac{d}{2}\rfloor-1$.
It follows that
$$Sd_k(T)=Sd_{k-3}(T)+l_{v_{k-2}}(S_{k-2})+l_{v_{k-1}}(S_{k-1})+l_{v_{k}}(S_k)\leq Sd_{k-3}(T)+3l_{v_{k-2}}(S_{k-2})$$
and
$Sd_{k-3}(T)\geq d+(|A|-2)\lfloor\frac{d}{2}\rfloor+(k-3-|A|)l_{v_{k-2}}(S_{k-2}).$
Therefore,
\begin{align*}
Sd_{k}(T)-\frac{k}{k-3}Sr_{k,3}(T)+\frac{6}{k-3}&\leq Sd_{k-3}(T)+3l_{v_{k-2}}(S_{k-2})-\frac{k}{k-3}Sr_{k,3}(T)+\frac{6}{k-3}\\
&= -\frac{3}{k-3}Sd_{k-3}(T)+3l_{v_{k-2}}(S_{k-2})+\frac{6}{k-3}\\
&\leq  -\frac{3}{k-3}(d+(|A|-2)\lfloor\frac{d}{2}\rfloor+(k-3-|A|)l_{v_{k-2}}(S_{k-2}))\\
&\quad +3l_{v_{k-2}}(S_{k-2}) +\frac{6}{k-3}\\
&= -\frac{3}{k-3}(d+(|A|-2)\lfloor\frac{d}{2}\rfloor)+\frac{3|A|}{k-3}l_{v_{k-2}}(S_{k-2}) +\frac{6}{k-3}\\
&\leq  -\frac{3}{k-3}(d+(|A|-2)\lfloor\frac{d}{2}\rfloor)+\frac{3|A|}{k-3}(\lfloor\frac{d}{2}\rfloor-1) +\frac{6}{k-3}\\
&\leq -\frac{3}{k-3}(d-2\lfloor\frac{d}{2}\rfloor)-\frac{3|A|}{k-3}+\frac{6}{k-3}\\
&\leq -\frac{3|A|}{k-3}+\frac{6}{k-3}\leq0.
\end{align*}
Thus, we obtain that $Sd_{k}(T)\leq\frac{k}{k-3}Sr_{k,3}(T)-\frac{6}{k-3}$.

\noindent{\bf Case 2. } $|A|=k-2$.

By $|A|=k-2$, we have $d_{T}(A)=Sd_{k-2}(T)$. By Lemma \ref{lem2}, we can choose $S_{k-1}$ such that $d_{T}(S_{k-1})=Sd_{k-1}(T)$ and $A\subset S_{k-1}$. If $|L(T)|=k-1$, then $d_{T}(S_{k-1})=Sd_{k-1}(T)=Sd_{k}(T)$. If $|L(T)|\geq k$, then by Lemma \ref{lem2}, there exists a $k$-subset $S_{k}$ of $V(T)$ satisfying $d_{T}(S_{k})=Sd_{k}(T)$ and $A\subset S_{k-1}\subset S_{k}$. Thus, we can choose $S$ such that $d_{T}(S)=Sd_k(T)$ and $A\subset S$. Then
%
\begin{align*}
Sd_k(T)=d_{T}(S)&\leq d_{T}(A)+\sum_{v\in S\setminus A}l_v(S)\\
&\leq d_{T}(A)+2(\lfloor\frac{d}{2}\rfloor-1)\\
&=d+(k-2)\lfloor\frac{d}{2}\rfloor-2.
\end{align*}
By $|A|= k-2$, we have
$Sr_{k,3}(T)=\varepsilon_{k,3}(S^{\prime \prime};T)=(k-4)\lfloor\frac{d}{2}\rfloor+\lceil\frac{d}{2}\rceil+1$. Then
\begin{align*}
Sd_k(T)-\frac{k}{k-3}Sr_{k,3}(T)+\frac{6}{k-3}&\leq d-2\lfloor\frac{d}{2}\rfloor+\frac{k}{k-3}(\lfloor\frac{d}{2}\rfloor-\lceil\frac{d}{2}\rceil)-\frac{k}{k-3}+\frac{6}{k-3}-2\\
&\leq-\frac{k}{k-3}+\frac{6}{k-3}-2\leq0.
\end{align*}
Thus, we obtain that $Sd_k(T)\leq\frac{k}{k-3}Sr_{k,3}(T)-\frac{6}{k-3}$.

\noindent{\bf Case 3. } $|A|\geq k-1$.

By Corollary \ref{cor2}, $Sd_k(T)\leq d+(k-2)\lfloor\frac{d}{2}\rfloor$.
By $|A|\geq k-1$, we have $Sr_{k,3}(T)=\varepsilon_{k,3}(S^{\prime \prime};T)=(k-3)\lfloor\frac{d}{2}\rfloor+2$. Then
\begin{align*}
Sd_k(T)-\frac{k}{k-3}Sr_{k,3}(T)+\frac{6}{k-3}\leq d-2\lfloor\frac{d}{2}\rfloor-2<0.
\end{align*}
Thus, we obtain that $Sd_k(T)<\frac{k}{k-3}Sr_{k,3}(T)-\frac{6}{k-3}$.

If $T\cong P_{3}(a,b;x)$, then $Sd_{k}(T)=kx+2$ and $Sr_{k,3}(T)=(k-3)x+2$. Thus $Sd_k(T)=\frac{k}{k-3}Sr_{k,3}(T)-\frac{6}{k-3}$.

Therefore, we conclude that $Sd_k(T)\leq\frac{k}{k-3}Sr_{k,3}(T)-\frac{6}{k-3}$, where the equality holds if $T\cong P_{3}(a,b;x)$. $\Box$

\section{Concluding Remarks}
~~~~In this paper, we obtain an upper bound on the Steiner $k$-diameter by its Steiner ($k,k^{\prime}$)-radius of a tree for $k\geq 3$, $k>k^{\prime}\geq1$, which generalizes the result for $k^{\prime}=1$ in \cite{diameter radius1}. Furthermore, we prove that $Sd_k(T)\leq \frac{k}{k-k^{\prime}}Sr_{k,k^{\prime}}(T)-\frac{k^{\prime}(k^{\prime}-1)}{k-k^{\prime}}$ for $k^{\prime}=2, 3$, and show that this bound is tight. Although it is only proved for $k^{\prime}=2,3$, we believe the inequality $Sd_k(T)\leq \frac{k}{k-k^{\prime}}Sr_{k,k^{\prime}}(T)-\frac{k^{\prime}(k^{\prime}-1)}{k-k^{\prime}}$ holds for all $k^{\prime}\geq 1$ and present it as a conjecture.

\begin{conjecture}
Let $n$, $k$, $k^{\prime}$ be integers with $n\geq k\geq 3$ and $k>k^{\prime}\geq 1$ and let $T$ be a tree with order $n$. Then
$$Sd_k(T)\leq \frac{k}{k-k^{\prime}}Sr_{k,k^{\prime}}(T)-\frac{k^{\prime}(k^{\prime}-1)}{k-k^{\prime}}.$$
\end{conjecture}


\end{sloppypar}
\end{document}